\documentclass{article}

\usepackage[T2A]{fontenc}
\usepackage[cp1251]{inputenc}
\usepackage[english,russian]{babel}
\usepackage[tbtags]{amsmath}
\usepackage{amsfonts,amssymb}

\usepackage{mathrsfs}

\usepackage{graphicx}

\numberwithin{equation}{section}

\newtheorem{theorem}{Теорема}
\newtheorem{proposition}{Предложение}
\newtheorem{lemma}{Лемма}

\newtheorem{definition}{Определение}

\newtheorem{proof}{Доказательство}
\newtheorem{remark}{Замечание}

\author{Д.В.~Гусев, И.А.~Иванов-Погодаев, А.Я.~Канель-Белов}

\title{Коллектив автоматов в конечно-порождённых группах}

\begin{document}


\maketitle



\begin{abstract}

Данная работа посвящена обхождению лабиринта коллективом конечных автоматов. Эта часть теории автоматов породила довольно широкий спектр различных задач \cite{kudr_kili_ush_1, andzhans_1}, в том числе связанных с задачами теории сложности вычислений и теорией вероятностей. Оказывается, что рассмотрение сложных алгебраических объектов, таких как бёрнсайдовы группы, может быть интересном в данном контексте. В работе будет показано, что граф Кэли конечно-порожденной группы нельзя обойти коллективом автоматов тогда и только тогда, когда она бесконечна и каждый ее элемент периодичен.

Библиография: 18 названий.

\end{abstract}


\footnotetext{Данная работа была проведена с помощью Российского Научного Фонда Грант N 17-11-01377.}

\section{Введение}
\label{sec1}

Есть довольно большое число вариаций задачи об обхождении автоматом лабиринта, но в целом она выглядит так: коллектив конечных автоматов двигается по рёбрам некоторого (возможно бесконечного) графа, необходимо выяснить смогут ли автоматы посетить все вершины графа.

Простейшим примером такой задачи, является обход $\mathbb{Z}$. Конечный автомат, двигающийся по решётке $\mathbb{Z}$ по некоторым внутренним правилам, не сможет обойти эту прямую, т.к. в какой-то момент зациклится. Однако можно рассмотреть коллектив из одного полноценного конечного автомата и двух автоматов-камней, которые не имеют внутреннего состояния и могут передвигаться только совместно с главным автоматом. Легко показать, что подобная система может обойти $\mathbb{Z}$ (главный автомат бегает между камнями и постепенно раздвигает их) \cite{andzhans_1}. В случае решёток $\mathbb{Z}^k, k > 1$ для обхода достаточно коллектива из автомата и трёх камней, причём с меньшим количеством камней обойти не получиться. Если из плоской решётки $\mathbb{Z}^2$ разрешить выкидывать некоторые вершины, то окажется, что автомата и трёх камней не достаточно для обхода \cite{kilibarda_3, kilibarda_4} таких лабиринтов, при этом пяти камней хватит \cite{szepietowski_1} (4 - открытый вопрос).

Возникает закономерный вопрос, существуют ли лабиринты, которые не обходятся подобными системами. Ответ на этот вопрос положительный. Например, можно построить бесконечную лабиринт-ловушку на решётке $\mathbb{Z}^3$ для любой системы автоматов \cite{kilibarda_ush_1}. Но построение таких лабиринтов обычно довольно громоздко. В работе мы предлагаем другой подход для построения подобных ловушек.

В качестве лабиринтов можно рассматривать графы Кэли конечнопорожденных групп. Такой подход даёт довольно интересные результаты. Оказывается, графы Кэли бесконечных свободных бёрнсайдовых групп, существование которых доказано Новиковым и Адаяном \cite{nov_adian_1}, нельзя обойти никакой системой конечных автоматов. Проблемы Бернсайдовского типа имеют большое значение в современной алгебре, поэтому кажется интересным, что такие серьёзные результаты находят применение в довольно удалённых от алгебры областях. В работе будет показано, что граф Кэли конечно-порожденной группы нельзя обойти коллективом конечных автоматов тогда и только тогда, когда она бесконечна и каждый ее элемент периодичен.

\section{Определения}
\label{sec2}

В этой секции мы введём формальные определения, которые нам понадобятся. Во-первых, мы определим понятие лабиринта. В нашем случае лабиринтом будет граф Кэли некоторой конечно-порождённой группы. Отметим, что возможны и другие определения лабиринта. Во-вторых, определим коллектив конечных автоматов и поведение коллектива в лабиринте. Также введём определение лабиринта"=ловушки.

Будем пользоваться основными понятиями и обозначениями теории из автоматов и графов, принятыми здесь \cite{kudr_podkol_ushmulich_1, kudr_alesh_podlol_1, harari_1}.

Пусть $G$ группа, $S = \lbrace s_1, s_2,...,s_n \rbrace$ -- множество образующих группы $G$, $S^{-1} = \lbrace s_1^{-1}, s_2^{-1}, ... , s_n^{-1} \rbrace$. Граф Кэли группы $G$ на образующих $S$ будем обозначать $\Gamma=\Gamma(G, S)$.  Для простоты будем считать, что $G$ - это вершины графа, также будет использовать групповые операции над вершинами, если это потребуется. За $E(\Gamma)$ обозначим множество рёбер графа. $\Gamma$ - будет лабиринтом в нашей задаче.

Под конечным автоматом $\mathfrak{A}$ будем понимать пятёрку $(A, Q, B,\phi, \psi)$, где $A, B, Q$ -- конечные множества (алфавиты): входной, выходной и алфавит состояний. $\phi : Q \times A \rightarrow Q$, $\psi : Q \times A \rightarrow B$ -- функции переходов и выходов. Если начальное состояние $q_0 \in Q$ фиксировано, то такой автомат $\mathfrak{A}_{q_0} = (A, Q, B,\phi, \psi, q_0)$ будем называть инициальным. Далее, если факт инициальности будет важен для нас, тогда будем обозначать автомат нижним индексом. Часто вместо конечного автомата будем говорить просто автомат.

Пусть задана группа $G$ с множеством образующих $S$.

\begin{definition}

Набор $\mathcal{A} = (\mathfrak{A}_1, \mathfrak{A}_2,...,\mathfrak{A}_m)$ называется коллективом $(G, S)$"=допустимых автоматов, где $\mathfrak{A}_i = (A_i, Q_i, B_i, \phi_i, \psi_i)$, для $i = 1,...m$, некоторые автоматы, такие что:

\begin{itemize}
\item $A_i = \lbrace a \in \prod\limits_{j=1}^m (\theta \cup Q_j) | \Pr_{i}(a) = \theta \rbrace $, где $\theta \notin \bigcup \limits_{i=1}^m Q_i$ -- некоторый выделенный элемент
\item $B_i = S \cup S^{-1} \cup \lbrace e \rbrace$
\end{itemize}

\end{definition}

Если у всех $m$ автоматов заданы начальные состояния $\vec{q_0} = (q_0^1, ... , q_0^m )$, то коллектив будем называть инициальным и обозначать $\mathcal{A}_{\vec{q}}$

Пусть $\Gamma(G, S)$ - граф Кэли. $\vec{v} = (v_1,...,v_m)$, $v_i \in G$ - его набор вершин; $\vec{q} = (q_1, ..., q_m)$, где $q_i \in Qb_i$ -- некоторое состояние автомата $\mathfrak{A}_i$. Тогда для для любого $i = 1,...,m$ определим функции:
\[
a_i(\vec{q}, \vec{v}) = (q_{i1}', q_{i2}', ... , q_{im}')
\]
где для $1 \leq j \leq m$, $q_{ij}' = q_j$, если $v_i = v_j$ и $i \neq j$, иначе $q_{ij} = \theta$.

\begin{definition}

Поведением коллектива $(G, S)$-допустимых автоматов \linebreak $\mathcal{A} = (\mathfrak{A}_1, \mathfrak{A}_2,...,\mathfrak{A}_m)$ с начальными состояниями $\vec{q_0} = (q_0^1, ... , q_0^m )$, $q_i \in Q_i$ в графе Кэли $\Gamma(G, S)$ с набором начальных вершин $\vec{v_0} = (v_0^1, ... , v_0^m)$, назовем последовательность \linebreak $\pi(\mathcal{A}_{\vec{q_0}}, \Gamma_{\vec{v_0}}) = (\vec{q_0}, \vec{v_0}), (\vec{q_1}, \vec{v_1}), ...,$ где $\vec{q_t} = (q_t^1, ... , q_t^m )$, $\vec{v_t} = (v_t^1, ... , v_t^m )$, такую, что

\begin{itemize}

\item $q_t^i \in Q_i$, $v_t^i \in G$, $(v_t^i, v_{t+1}^i) \in E(\Gamma)$ или $v_t^i = v_{t+1}^i$
\item $\phi_i(q_t^i, a_i(\vec{v_t}, \vec{q_t})) = q_{t+1}^i$

\item $\psi_i(q_t^i, a_i(\vec{v_t}, \vec{q_t})) = g$, где $v_t^i g = v_{t+1}^i, g \in S \cup S^{-1} \cup e$

\end{itemize}

\end{definition}

То есть, у нас есть коллектив автоматов, который перемещается по лабиринту. Они могут взаимодействовать друг с другом, когда находятся в одной вершине. Формально - следующее состояние и направление хода одного автомата, зависит от состояний автоматов, которые находятся с ним в одной вершине.


Пусть,

\[
 \mathrm{Int}_i (\mathcal{A}_{\vec{q_0}}, \Gamma_{\vec{v_0}}) = \bigcup\limits_{t=0}^{\infty} \lbrace v_t^i \rbrace, \forall i = 1,...,m
\]

\[
 \mathrm{Int} (\mathcal{A}_{\vec{q_0}}, \Gamma_{\vec{v_0}}) = \bigcup\limits_{t=1}^{m} \mathrm{Int}_i (\mathcal{A}_{\vec{q_0}}, \Gamma_{\vec{v_0}})
\]

\[
\mathrm{Fr} (\mathcal{A}_{\vec{q_0}}, \Gamma_{\vec{v_0}}) = V \setminus \mathrm{Int} (\mathcal{A}_{\vec{q_0}}, \Gamma_{\vec{v_0}})
\]

Будем говорить, что $\mathcal{A}_{q_0}$ обходит $\Gamma_{v_0}$, если $\mathrm{Fr} (\mathfrak{A}_{q_0}, \Gamma_{\vec{v_0}}) = \emptyset$. Говорим, что $\mathcal{A}_{q_0}$ сильно обходит $\Gamma$, если он обходит его, начиная из любого набора $\vec{v_0}$ вершин лабиринта. $\Gamma_{v_0}$ -- ловушка для $\mathcal{A}_{q_0}$, если автомат не обходит граф; $\Gamma$ -- сильная ловушка для $\mathcal{A}_{q_0}$, если автомат не может обойти лабиринт из любого набора начальных вершин.

Иногда выделяют особый вид коллектива - один основной автомат и $k$ автоматов"=камней. В данном виде коллектива автоматы-камни могут двигаться только совместно с основным автоматом, фактически основной автомат их двигает по лабиринту, а сами по себе они перемещаться не могут.

\section{Построение сильной ловушки для любого коллектива автоматов}

В этой секции мы построим сильную ловушку для любого коллектива автоматов. Ловушка будет выглядеть как граф Кэли некоторой группы. Сама группа будет и, соответственно, лабиринт бесконечны, но любой коллектив автоматов обойдёт лишь конечную область в нём, что мы и докажем.

\subsection{Проблема Бёрнсайда и её решения}

Проблема Бёрнсайда о периодических группах фиксированного периода была поставлена Бёрнсайдом в 1902 году в следующей форме \cite{burnside_1}

\begin{proposition}

Пусть группа $G$ имеет $m$ независимых порождающих элементов \linebreak $a_1, a_2, ..., a_m$ и для любого элемента $x \in G$ выполнено соотношение $x^n=1$, где n — данное целое число. Будет ли определенная таким образом группа конечной, и если да, то каков ее порядок?

\end{proposition}

Сейчас группы, определенные $m$ порождающими и соотношением $x^n=1$, называют свободными бёрнсайдовыми группами ранга $m$ и периода $n$ (экспоненты $n$). Обычно они обозначаются как $B(m,n)$.

Понятно, что вопрос нетривиален для случая $m > 1$. Самим Бёрнсайдом показана конечность $B(m,n)$ для $n \leq 3$ и любого $m$ и $B(2,4)$ \cite{burnside_1}. Сановым \cite{sanov_1} показана конечность $B(m,n)$ для $n=4$ и любого, М.Холлом \cite{hall_1} для $n=6$ и любого.

В 1964 году Голод и Шафаревич \cite{golod_1} доказали, что существуют бесконечные $2$ - порожденные периодические группы с неограниченными периодами элементов. В 1968 году Новиков и Адян представили отрицательное решение проблемы Бёрнсайда \cite{nov_adian_1} для любого нечётного периода $n > 4381$ и любого $m > 1$. В \cite{adian_1} Адяном решение упрощено и доказано, что $B(m,n)$ бесконечны для любого нечётного $n > 665$ и любого $m > 1$. Недавно Адяном было объявлено о улучшении оценки до $n \ge 101$ \cite{Adyan1}.

В 1994 году Иванов \cite{ivanov_1} доказал существование бесконечных Бернсайдовых групп для случая четного периода при $n \geq 2^{48}$. В 1996 Лысёнок \cite{Lys96} улучшил результат для четных периодов до $n \geq 8000$.

Таким образом, мы имеет бесконечные группы, порождённые конечным количеством элементов, при этом период каждого элемента равномерно ограничен.

\subsection{Определение лабиринта}

Пусть $G$ бесконечная конечно-порождённая\-группа, у которой нет элементов бесконечного порядка. $S = \lbrace s_1, s_2,...,s_n \rbrace$ -- множество образующих группы $G$, среди которых нет повторяющихся и обратных. $\Gamma(G, S)$ - граф Кэли для этой группы и образующих.

Пусть $g = t_1 t_2 ... t_n$, где $g \in G$ и $t_i \in S \cup S^{-1} $, причем $n$ минимальное число с такими свойствами, т.е $g$ нельзя представить меньшим количеством элементов из $S \cup S^{-1}$. Тогда будем обозначать $d(g) = n$. Заметим, что $d(g) = d(g^{-1})$. Также, легко видеть, что кратчайшее расстояние в лабиринте между двумя вершинами $v_1$ и $v_2$, равно $d(v_1^{-1} v_2) = d(v_2 v_1^{-1})$. $M(r)$ - общий период всех элементов $g$, для которых $d(g) \leq r$. В силу отсутствия элементов бесконечного порядка $M(r)$ конечно и определено для всех $r$.

\subsection{Единичный автомат}

Рассмотрим сначала случай с одним автоматом. Покажем, что единичный автомат сможет обойти только ограниченную часть $\Gamma$.

Пусть $\mathfrak{A}$ -- $(G, S)$-допустимый коллектив из одного автомата. Заметим, что ход автомата в конкретной вершине определяется исключительно его состоянием, т.к. других автоматов нет. Пусть этот автомат $\mathfrak{A}$ начинает движение из некоторой вершины $v_0$. Тогда, его движение будет довольно простым.

\begin{lemma}

Поведение автомата $\mathfrak{A}_{q_0}$ в $\Gamma$ c $|Q|$ состояниями будет обладать следующими свойствами:

\begin{itemize}

\item На начальной стадии автомат сделает $U < |Q|$ ходов, с неповторяющимися состояниями;
\item Далее состояния будут повторяться с периодом $T \leq |Q|$;
\item Каждые $M(T)$ периодов вершины, посещаемые автоматом будут повторятся;

\end{itemize}

\end{lemma}

\begin{proof}

Так как количество состояний конечно, какое-то состояние автомата повториться как минимум два раза. Пусть первое такое состояние $q_1$ и между первым и вторым появлением сделано $T$ ходов (очевидно $T \leq |Q|$), а перед $q_1$ было $U < |Q|$ ходов. Заметим, что текущее состояние однозначно определяет следующее, таким образом после второго появления $q_1$ последовательность состояний будет такой же, как и после первого появления. Отсюда получаем зацикленность состояний с периодом $T$. Обозначим через $v_1$ вершину лабиринта, в которой автомат впервые оказался в состоянии $q_1$.

Пусть $s_1', s_2', ... , s_T'$ -- направления движения автомата в цикле, $s_i' \in S \cup S^{-1}$. Пусть $g_T = s_1' s_2' ... s_T'$, $g_T \in G$. Тогда через $M(T)$ циклов после посещения вершины $v_1$ автомат окажется в вершине $v_1 g_T^M = v_1$, то есть вернётся обратно. Таким образом, сначала автомат посетит не более $U$ вершин, а дальше в цикле будет посещать не более $TM(T)$ вершин.

\end{proof}

Из этой леммы легко видеть, что $\Gamma$ -- сильная ловушка для любого одиночного автомата.

\subsection{Коллектив}

Теперь рассмотрим случай коллектива автоматов. Пусть \linebreak $\mathcal{A}_{\vec{q_0}} = (\mathfrak{A}_1, \mathfrak{A}_2,...,\mathfrak{A}_m)$ -- $(G, S)$-допустимый коллектив автоматов. Основная наша цель -- доказать, что построенный нами лабиринт -- сильная ловушка для любого такого коллектива, причём коллектив сможет обойти только конечное число вершин лабиринта.

Пусть $|Q_1|, ... |Q_m|$ -- количество состояний у автоматов. Обозначим через $Q_A$ максимум из $|Q_i|$ $i = 1,...,m$.

Назовём состоянием коллектива $\mathcal{A}$ в момент времени $t$ набор $I = (q_t^1, q_t^2, ... q_t^m, F_t) $, где $q_t^i \in Q_i$ -- состояние конкретного автомата в момент $t$; $F_t = \lbrace F_t^1, F_t^2, ... , F_t^k \rbrace$, $k \leq m$, где $F_t^i$ -- множество номеров всех автоматов находящихся в какой-то вершине в данный момент. Причём $\bigcup \limits_{i=1}^k F_t^i = \lbrace 1,...,m \rbrace$ и один индекс принадлежит только одному $F_t^i$. Если все автоматы находятся в разных вершинах, то $F_t^i$ состоят их одного индекса каждый и их $m$ штук. То есть, состояние коллектива -- это состояния всех автоматов и разбивка автоматов на группы стоящих в одной вершине. Далее факт нахождения в одной вершине двух или более автоматов будем называть встречей.

Положением коллектива в лабиринте в момент времени $t$, назовём набор вершин
\linebreak $(v_t^1, v_t^2, ... ,v_t^m)$, в котором находятся автоматы. Заметим, что пара положение"=состояние однозначно все определяет последующее поведение коллектива автоматов.

\begin{theorem}

Поведение коллектива автоматов $\mathcal{A}_{\vec{q_0}} = (\mathfrak{A}_1, \mathfrak{A}_2,...,\mathfrak{A}_m)$ в $\Gamma$ в вектором начальных вершин $\vec{v_0}$ обладает следующим свойством: существуют такие $O_m$ и $H_m$ зависящее только от $m$, $Q_A$ и функции $M(r)$, что:

\begin{itemize}

\item Состояния коллектива $\mathcal{A}$ любые в $H_m$ подряд идущих моментов времени однозначно определяют следующее состояние.
\item На начальной стадии коллектив автоматов сделает  $U \leq O_m$ ходов, после которых состояния будут повторяться с периодом $T \leq O_m$;
\item Каждые $M(T)$ периодов, положение коллектива будет повторяться, т.е. после первых $U$ ходов пара положение-состояние коллектива $\mathcal{A}$ будет повторяться с периодом $M(T)T$.

\end{itemize}

\end{theorem}

\begin{proof}

Будем доказывать индукции по числу автоматов. В случае с одним автоматом, состояние коллектива фактически состояние автомата и $O_1 = Q_A$, $H_1 = 1$, что было доказано ранее. Таким образом база индукции есть.

Пусть $m=l$, тогда по предположению индукции для всех $m < l$ утверждение верно. $H_l$ и $O_l$ будем выражать через $l$, $M$, $Q_A$, $H_i$ и $O_i$ для $i < l$. Поскольку $H_i$ и $O_i$ по предположению зависят только от $i$, $Q_A$ и функции $M(r)$, $H_l$ и $O_l$ тоже можно будет выразить через $l$, $Q_A$ и функцию $M(r)$.

\begin{lemma}
Пусть существуют две группы автоматов размерами $a$ и $b$, $a + b = l$. Пусть $O_a$ и $O_b$ константы из индукции. Тогда, если никакие два автомата из этих групп (один из первой, один из второй) не встречались $h = \mathrm{max} (O_a, O_b) + M(O_a) M(O_b) O_a O_b + 1$ ходов подряд, то и в дальнейшем никакие два не встретятся.
\end{lemma}

\begin{proof}
Заметим, что эти $h$ ходов, группы не влияют на движения друг друга, поэтому эти две группы мы можем в течении $h$ ходов рассматривать как два коллектива в лабиринте $L$ c некоторыми начальными состояниями и вектором начальных вершин. Заметим, так как $a,b<l$, то к этим коллективам применимы предположения индукции (но только на $h$ первых ходов, в дальнейшем в худшем случае автоматы из разных групп встретятся и независимость пропадёт).

Тогда, после $\mathrm{max} (O_a, O_b)$ ходов (могут и раньше), состояния автоматов в обоих группах будут повторяться циклами длинной $T_a \leq O_a$ и $T_b \leq O_b$ и каждые $M(T_a)$ и $M(T_b)$ таких циклов положения групп автоматов повторяются.

Заметим, что за следующие $M(T_a) M(T_b) T_a T_b + 1$ ходов первая группа будет находиться в не более $M(T_a) T_a$ различных парах положение-состояние, аналогично вторая в не более $M(T_b) T_b$ различных парах. Тогда существуют два момента времени $t_1, t_2$ из этих $M(T_a) M(T_b) T_a T_b + 1$ ходов, когда положение первого коллектива -- $\vec{v_a}$, состояние -- $(\vec{q_a}, F_a)$, положение второго коллектива -- $\vec{v_b}$, состояние -- $(\vec{q_b}, F_b)$. Заметим, что из эти двух пар положение-состояние однозначно получается пара положение-состояние всего коллектива. Таким образом, пара положение-состояние всего коллектива из $m$ автоматов повторилась, то есть оно и дальше будет повторяться с периодом $T_{ab} = t_2 - t_1 \leq M(T_a) M(T_b) T_a T_b$. Поскольку, между $t_1$ и $t_2$ группы автоматов не встречались, они и дальше не будут встречаться. Заметим, что в силу монотонности $M(r)$ $M(T_a) M(T_b) T_a T_b \leq M(O_a) M(O_b) O_a O_b$, откуда следует требуемое.

\end{proof}

Докажем следующую лемму, являющуюся подпунктом теоремы.

\begin{lemma}

Положим $H_l$ минимальным натуральным числом со следующими свойствами:
\[
H_l \geq \mathrm{max} (O_i, O_j) + M(O_a) M(O_b) O_i O_j + 1, \forall i > 0, j > 0, i + j = l
\]
\[
	H_l \geq  H_i, \forall i < l
\]

Тогда, по последовательности $I_1, I_2, ..., I_{H_l}$ из $H_l$ состояний коллектива можно однозначно определить состояние $I_{H_l+1}$.

\end{lemma}

\begin{proof}

Заметим, что по этим $H_l$ состояниям мы можем определить, есть ли две группы автоматов, автоматы из которых не встречаются друг с другом за эти $H_l$ ходов (в состояниях записана информация о встречах автоматов).

Пусть есть две такие группы $A$ и $B$ c размерами $a$ и $b$. Тогда в силу выбора $H_l$ эти группы удовлетворяют условию предыдущей леммы. Следовательно, автоматы из этих групп не встретятся после этих $H_l$ ходов, а, следовательно, и на ходу $H_l+1$. Заметим, что состояние коллектива из $l$ автоматов однозначно определяется состоянием коллектива $A$ и коллектива $B$ и наоборот (если их рассматривать как коллективы автоматов в $L$) на этих $H_l + 1$ ходу.

Тогда, $I_{H_l+1}$ -- состояние коллектива $\mathcal{A}$ из $l$ автоматов, однозначно определяется состоянием коллектива $A$ и коллектива $B$ на $H_l+1$-ом ходу. В свою очередь, состояния коллективов $A$ и $B$ на этом ходу, однозначно определяются предыдущими своими $H_a$ и $H_b$ состояниями, а значит и предыдущими $H_l$ состояниями. Следовательно $I_{H_l+1}$ однозначно определяется через предыдущие $H_l$ состояний коллективов $A$ и $B$, а значит и $\mathcal{A}$, что и требовалось.

Пусть таких двух групп нет. Далее все моменты времени лежат в интервале от $1$ до $H_l$ (если не оговорено обратное). Заметим, что по состояниям $I_1, I_2, ..., I_{H_l}$ можно определить, направление движения автомата $\mathfrak{A}_i$ в момент $t$ (будем обозначать его $g_t^i \in D(S)$). Следовательно, если известная вершина $v_{t_1}$, в которой $\mathfrak{A}_i$ находился в момент $t_1$, то и известна вершина $v_{t_2}$, в которой он находился любой другой момент $t_2$, причём $v_{t_2} = v_{t_1}g$, где $g \in G$ произведение каких-то $g_j^i$, либо их обратных.

Пусть $\mathfrak{A}_1$ в момент времени $1$ находится в вершине $v_0$. Тогда в остальные моменты $2 \leq t \leq H_l$ он находится в вершинах $v_0 h_t^1$, где $h_t^1 \in G$ однозначно определяется через $I_1, I_2, ..., I_{H_l}$. В силу предположения о том, что разделённых групп не существует, есть автомат, для удобства $\mathfrak{A}_2$, с которым он встретится в момент $t_1$, в вершине $v_0 h_{t_1}^1$. Тогда, в момент $t_1$ известно положение $\mathfrak{A}_2$, а значит известны вершины в которой $\mathfrak{A}_2$ находится во все остальные моменты времени $t$, причём они имеют вид $v_0 h_{t}^2$, где $h_{t}^2$ однозначно выражается через $I_1, I_2, ..., I_{H_l}$ ($h_{t_1}^1 = h_{t_1}^2$, а остальные получаются из $h_{t_1}^2$ с помощью $g_j^2$).

Аналогично, найдем автомат $\mathfrak{A}_3$, с которым встречается хотя бы один из двух первых, для него положения вершин будут определятся аналогично. Так будем продолжать и далее. Заметим, что в силу условия о не существовании двух не встречающихся групп, мы всегда найдем новый автомат, с которым есть встречи у уже рассмотренных. Таким образом, любой $\mathfrak{A}_i$ из $m$ автоматов в момент $t$ находится в вершине $v_0 h_t^i$, причём $h_{t}^i$ однозначно выражается через $I_1, I_2, ..., I_{H_l}$.

Рассмотрим момент времени $H_l+1$. Состояния каждого из $m$ автоматов, однозначно получаются из $I_{H_l}$. Пусть в момент $H_l+1$ $\mathfrak{A}_i$ находится в вершине $v_0 h_{H_l}^i g_i$, $g_i \in G$. Заметим, что для определения, стоят ли $\mathfrak{A}_i$ и $\mathfrak{A}_j$ достаточно сравнить на равенство $h_{H_l}^i g_i$ и $h_{H_l}^j g_j$. $h_{H_l}^i$ и $h_{H_l}^j$ однозначно определяется через $I_1, I_2, ..., I_{H_l}$, $g_i$ и $g_j$ через $I_{H_l}$. Следовательно, разбиение на группы стоящие в одной вершине на шаге $H_l + 1$ однозначно определяется через $I_1, I_2, ..., I_{H_l}$, а значит и $I_{H_l+1}$ однозначно через них определяется. Что нам и требовалось. Итого пункт про $H_l$ выполнен.

\end{proof}

Посчитаем, сколько различных состояний $I = (\vec{q}, F)$ может быть у коллектива $\mathcal{A}$. Различных $\vec{q}$ не более $Q_A^l$ (т.к. состояний у каждого из автоматов не более $Q_A$). Количество различных $F$  грубо оценим сверху числом $l^l$. Действительно, разложить $l$ индексов автоматов по $l$ множествам можно не более чем $l^l$ способами, и каждому такому способу соответствует только одно $F$ (разбиение автоматов по группам, стоящим в одной вершине). Итого, различных состояний коллектива не более чем $(Q_A l)^l$.

Тогда заметим, что существует не более $(Q_A l)^{lH_l}$ различных блоков из $H_l$ подряд идущих состояний коллектива. Положим $O_l = (Q_A l)^{lH_l} + H_l$, рассмотрим первые $O_l$ ходов коллектива $\mathcal{A}$ в $L$. В этих ходах можно выделить $(Q_A l)^{lH_l} + 1$ блок из $H_l$ подряд идущих состояний коллектива. Тогда, два блока начинающиеся в ходы $t_1$ и $t_2$ будут одинаковыми. Так как каждый такой блок однозначно определяет последующую последовательность состояний коллектива, мы получаем, что, начиная с момента $t_1$, состояния коллектива зацикливаются с периодом $t_2 - t_1$. Таким образом, $U$ из условия теоремы равно $t_1 \leq O_l$ (если первый ход считается нулевым), а $T = t_2 - t_1 \leq O_l$. Получаем, что второй пункт теоремы выполнен, нужное $O_l$ найдено

Для доказательства теоремы осталось показать, что каждые $M(T)$ получившихся циклов, положения коллектива будут повторяться. Действительно, пусть любой  $\mathfrak{A}_i$ на каком-то ходе $t \geq t_1$ находится в вершине $v_t$. Пусть за каждые $T$ ходов он сдвигается на $g_i \in G$. $g_i$ слово длины не более чем $T$, поэтому через $M(T)$ ходов $\mathfrak{A}_i$ попадёт в вершину $v_t g_i^M(T) = v_t$, что нам и требуется. Таким образом после первых $U$ ходов пара положение-состояние всего коллектива $\mathcal{A}$ будет повторяться с периодом $M(T)T \leq M(O_m) O_m$.

\end{proof}

Таким образом, мы показали, что любой коллектив автоматов стартующий из любых вершин обойдёт только конечную часть $\Gamma$.

\section{Об обходе графов Кэли конечно-порожденных групп}

Можно сформулировать следующее утверждение об обходе коллективом автоматов графов Кэли конечно-порожденных групп.

\begin{theorem}
Граф Кэли конечно-порожденной группы нельзя обойти коллективом автоматов тогда и только тогда, когда она бесконечна и в ней нет элемента бесконечного порядка.
\end{theorem}

Невозможность обхода бесконечной конечно-порождённой периодической группы мы показали выше. Приведём здесь набросок доказательства существования обхода групп с элементами бесконечного порядка.

\begin{lemma}
Пусть $G$~-- конечно порожденная группа, содержащая элемент $g\in G$ бесконечного порядка. Тогда ее граф Кэли $\Gamma$ может обойти коллектив из одного автомата и трёх автоматов-камней.
\end{lemma}

\begin{remark}
Минимально возможное число состояний автомата зависит от элемента $g$ и графа $\Gamma$.
\end{remark}

\begin{proof}

Стандартно и полностью аналогично теореме Аджанса \cite{andzhans_1} об обходе бесконечного дерева автоматом с тремя камнями, мы приведем его набросок. Пусть $s$~-- число образующих группы $G$, тогда $\Gamma$ есть естественная проекция $\pi(T_{2s})$ свободного дерева $T_{2s}$ степени $2s$ (удвоение происходит из-за наличия обратных элементов к образующим группы). Расположим три камня вдоль степеней элемента $g$. Пусть они расположены в точках $h, g^kh, g^{k+l}h$, автоматом между ними. Тогда с величинами $k,l$ автоматом может оперировать как со счетчиками, а два счетчика задают машину Минского \cite{andzhans_1,kudr_podkol_ushmulich_1}, для которой можно написать программу обхода $T_{2s}$. Тогда произойдет и обход $\pi(T_{2s})$.

Остается пояснить, как осуществить перемещение системы автомат -- камни вдоль ребра $v$ графа Кэли. Рассмотрим момент, когда в одном из счетчиков стоит $0$, т.е. камни $K_1$ и $K_2$ стоят в одной клетке $h$, а камень $K_3$~-- в клетке $g^mh$.

Далее автомат делает следующие шаги:

1) Переносит $K_1$ вдоль $v$, затем переходит обратно к $K_2$ (сдвигается на $v^{-1}$);

2) Идет вдоль $g$, пока не дойдет до $K_3$, переносит $K_3$ на $g^{-1}$, далее идет назад вдоль $g$ без $K_3$, пока не дойдет до  $K_2$;

3) Переходит по $v$, далее идет по $g$ пока не дойдет до $K_1$, переносит $K_1$ вперед на $g$;

4) Далее делает шаг по $g$ назад, потом проверяет проверяет соседние клетки на наличие $K_2$, повторяет эту операцию, пока не найдет $K_2$;

5) Повторяет шаги 2-4,  пока $K_3$ не окажется в одной клетке с $K_2$. После этого переносит $K_3$ и $K_2$ сдвигает на $v$ что завершает сдвиг системы.

6) Далее запускает машину Минского чтобы учесть информацию о сдвиге на $v$, и вычислить следующий элемент сдвига.

Таким образом, имея в распоряжении машину Минского с записанной в ней программой посещения графа, и перетаскивая ее за собой, автомат обходит весь граф.


\end{proof}


\end{document}